\documentclass[11pt]{amsart}
\usepackage{amsmath}
\usepackage{amsfonts}
\usepackage{amssymb}
\usepackage{amsthm}
\usepackage{marginnote}
\usepackage{graphicx}
\usepackage[section]{placeins}

\newtheorem{thm}{Theorem}

\newtheorem{cor}[thm]{Corollary}

\theoremstyle{definition}

\newtheorem{rem}[thm]{Remark}

\newcommand{\Z}{\mathbb{Z}}
\newcommand{\C}{\mathbb{C}}

\newcommand{\R}{\mathbb{R}}
\renewcommand{\S}{\mathrm{S}}
\renewcommand{\epsilon}{\varepsilon}
\newcommand{\del}{\partial}
\renewcommand{\d}{d}

\DeclareMathOperator{\can}{can}

\title[Lagrangian fillings and complicated\dots]
{Lagrangian fillings and complicated Legendrian unknots}
\author{Sylvain Courte}
\address{Institut Fourier, Universit\'e Grenoble Alpes et CNRS, 100,
rue des maths, 38610, Gi\`eres, France}
\author{Tobias Ekholm}
\address{Department of mathematics, Uppsala University, 751 06, Sweden\newline 
\indent Institut Mittag-Leffler, Aurav 17, 182 60 Djursholm, Sweden}
\thanks{TE was partially supported by the Knut and Alice Wallenberg
Foundation and the Swedish Research Council. SC was partially supported by
the ANR project MICROLOCAL (ANR-15CE40-0007-01)}

\begin{document}

\maketitle

\begin{abstract}
An exact Lagrangian submanifold $L$ in the symplectization of standard
contact $(2n-1)$-space with Legendrian boundary $\Sigma$ can be glued
to itself along $\Sigma$. This gives a Legendrian embedding $\Lambda(L,L)$
of the double of $L$ into contact $(2n+1)$-space. We show that the Legendrian
isotopy class of $\Lambda(L,L)$ is determined by formal data: the manifold $L$
together with a trivialization of its complexified tangent bundle. In particular,
if $L$ is a disk then $\Lambda(L,L)$ is the Legendrian unknot.
\end{abstract}

\section{Introduction}\label{sec:intr}
If $\Sigma$ is a closed Legendrian $(n-1)$-submanifold in standard contact
$\R^{2n-1}$ and $L_{0}$ and $L_{1}$ are exact Lagrangian fillings of $\Sigma$
in the symplectization $\S\R^{2n-1}=\R\times\R^{2n-1}$, then one can glue
compact versions of $L_{0}$ and $L_{1}$ along $\Sigma$ to form a Legendrian
submanifold $\Lambda(L_{0},L_{1})$ in standard contact $\R^{2n+1}$. The
Legendrian $\Lambda(L_{0},L_{1})$ has the same Reeb chords as the original
Legendrian $\Sigma$, but with Conley-Zehnder indices shifted up by $1$.
In \cite{ekholm_2016}, the second author used this construction for Legendrian spheres
$\Sigma$ with distinct disk fillings $D_{0}$ and $D_{1}$ and originally claimed that $\Lambda(D_0,D_1)$ were non-loose. These spheres were however later shown by Murphy to be loose. We present her proof in Section \ref{sec:loose}. 

The Chekanov-Eliashberg algebra of a loose Legendrian is trivial. In \cite{ekholm_2016}, it was shown that the Chekanov-Eliashberg dg-algebra for the double
$\Lambda(D_{0},D_{0})$ obtained by gluing a Lagrangian disk filling to itself
is quasi-isomorphic to that of the standard Legendrian unknot, but the analysis there left open the question of whether
$\Lambda(D_0,D_0)$ is actually Legendrian isotopic to the unknot. We prove
here that this is indeed the case, see Corollary \ref{cor:unknot}.
We point out that the Reeb chord configuration of $\Sigma$ and hence of the double $\Lambda(D_0,D_0)$ can be
involved and that when this is the case, the Legendrian isotopy to the unknot
must be rather complicated.
 
Disk fillings are special cases of Lagrangian fillings of more general topology
where a similar result holds that we describe next. To state the result, recall that a
Lagrangian embedding of a manifold $L$ in $\S \R^{2n-1}$
induces a trivialization of $T L \otimes \C$ that is well-defined up to homotopy.

\begin{thm}\label{t:main}
If $L\subset \S\R^{2n-1}$ is an embedded Lagrangian submanifold with Legendrian boundary,
then the Legendrian isotopy class of $\Lambda(L,L) \subset \R^{2n+1}$ is
determined by the induced trivialization of $TL\otimes\C$.
\end{thm}

\begin{cor}\label{cor:unknot}
If $D\subset \S \R^{2n-1}$ is an embedded Lagrangian \emph{disk} with Legendrian
boundary, then $\Lambda(D,D)$ is the Legendrian unknot.
\end{cor}

\subsection*{Acknowledgements} We thank Yang Huang for telling us how to simplify and strengthen our original version of Corollary 2. We also thank Emmy Murphy for proving the looseness results in Section \ref{sec:loose} and for allowing us  to present these arguments here.

\section{Gluing two Lagrangian fillings}

In this section we review the gluing construction from \cite[Section 2.1]{ekholm_2016}
in notation convenient for our study here. 

\subsection{Basic notation for contact structures and exact Lagrangians}\label{sec:basics}
Let $n\geq 1$ and consider $\R^{2n-1}$ with coordinates
$(q_{1}, p_{1},\dots, q_{n-1},p_{n-1},\zeta)$ and the contact structure 
\[\eta=\ker\left(\d\zeta - \sum_{i=1}^{n-1}p_{i} \d q_i\right),\]
and its symplectization $\S \R^{2n-1}=\R \times \R^{2n-1}$ with the canonical Liouville form
\[\lambda_{\can}=e^t\left(\d\zeta - \sum_{i=1}^{n-1}p_{i} \d q_i\right),\]
where $t$ is a standard coordinate on the additional $\R$-factor in $\S\R^{2n-1}$.
We write $X_{\can}$ for the canonical Liouville vector field on $\S \R^{2n-1}$, $X_{\can}=\del_{t}$,
and $\omega_{\can}$ for the canonical symplectic form, $\omega_{\can}=\d \lambda_{\can}$.

Recall that an \emph{exact Lagrangian cobordism} in $\S \R^{2n-1}$ is a proper embedding
$\phi\colon L \to \S \R^{2n-1}$ such that $\phi^*\lambda_{\can}=\d h$ for some compactly
supported function $h\colon L \to \R$ ($h$ is called a primitive of $\lambda_{\can}$).
This implies that for all sufficiently large $a>0$, $\phi(L)$ coincides with the
symplectization of a closed Legendrian submanifold $\Sigma_-$ (resp.~$\Sigma_+$)
$(-\infty,-a)\times\Sigma_{-}$ in $(-\infty,-a)\times \R^{2n-1}$
(resp. $(a,\infty)\times\Sigma_{+}$ in $(a,+\infty) \times \R^{2n-1}$). (In the case when $\Sigma_-$ is empty, we say that
$L$ is an exact Lagrangian \emph{filling} of $\Sigma_+$.) It follows in
particular that $L$ is diffeomorphic to a compact cobordism $\bar L$ with cylindrical ends
$(-\infty,-a]\times \Sigma_-$ and $[a,+\infty) \times \Sigma_+$ attached to the boundary
$\partial \bar L=\Sigma_{-}\cup\Sigma_{+}$. 

\begin{rem} 
We will also use the notion \emph{exact Lagrangian cobordism with Legendrian boundary}
when the Lagrangian $L$ and the ambient exact symplectic manifold $W$ are compact manifolds with boundary. This then means that $L$ is tangent to the Liouville
vector field of $W$ near the boundary $\del W$. Let $L\subset W$ be a compact Lagrangian cobordism with Legendrian boundary consisting of $\Lambda_{+}\subset\partial_{+}W$ and $\Lambda_{-}\subset\partial_{-}W$, where $\partial_{+}W$ (resp.~ $\partial_{-}W$) denotes the convex (resp.~concave) boundary of $W$ where the Liouville vector fields points out of (resp.~into) $W$. Then $(W,L)$ can be canonically \emph{completed}
by attaching cylindrical ends $\left([0,\infty)\times \partial_{+}W,[0,\infty)\times \Lambda_{+}\right)$ (resp.~$\left((-\infty,0]\times \partial_{-}W,(-\infty,0]\times \Lambda_{-}\right)$) to the convex (resp.~concave) boundary of $(W,L)$. In this way, the compact and non-compact cases are canonically related and we will use the same notation in both cases. For example, a Hamiltonian isotopy between compact Lagrangian
cobordisms, means a Hamiltonian isotopy of contact type at infinity between their completions. 
\end{rem}

Consider next $\R^{2n+1}$ with coordinates
$(x_1,y_1,\dots,x_{n},y_{n},z)$ and introduce polar coordinates on the last coordinate pair: 
\[(x_n,y_n)=(r\cos\theta,r\sin \theta),\quad (r,\theta)\in [0,\infty)\times \R/2\pi\Z.\]
Using these coordinates, we equip $\R^{2n+1}$ with the contact form
\[\alpha=\d z -\sum_{i=1}^{n-1} y_i\d x_i + \frac{1}{2} r^2 \d \theta,\]
and the corresponding contact structure $\xi = \ker \alpha$. 

We will occasionally discuss front projections, and for that purpose we
also consider the contact form 
\[\alpha'=\d z - \sum_{i=1}^{n-1} y_i\d x_i - y_n \d x_n,\]
with corresponding contact structure $\xi'=\ker \alpha'$.
The diffeomorphism $\psi\colon \R^{2n+1} \to \R^{2n+1}$,
\[\psi(x_1,y_1,\dots, x_n,y_n,z) = \left(x_1,y_1,\dots,x_n,y_n,z+\frac{x_ny_n}{2}\right),\]
satisfies $\psi^*\alpha'=\alpha$ and is thus a contactomorphism intertwining the
contact structures $\xi$ and $\xi'$.  

\subsection{The gluing construction}\label{sec:glue}
We now recall the construction from \cite[Section 2.1]{ekholm_2016} (with
slight notational changes). Pick a smooth function $\rho\colon [-1,1] \to \R$ with
the following properties:
\begin{enumerate}
\item $\rho$ is even,
\item $\rho(u)=1-|u|$ near $u=\pm 1$,
\item $\rho >0$ and its derivative $\rho'$ satisfies $\rho'<0$ on $(0,1)$.
\end{enumerate}
Use coordinates on $\R^{2n+1}$ as in Section \ref{sec:basics} and consider
the hypersurface $W_{\rho}\subset\R^{2n+1}\cap\{0<x_{n}<1\}$ given by,
\[W_\rho=\left\{z=\left(\frac{\rho(x_n)}{\rho'(x_n)} -\frac{x_n}{2}\right) y_n\right\}.\]
The restriction of $\alpha$ to $W_\rho$ is then
\[\lambda_\rho:=\alpha|_{W_{\rho}} = \frac{\rho(x_n)}{\rho'(x_n)} \d y_n -
y_n \frac{\rho(x_n)\rho''(x_n)}{(\rho'(x_n))^2} \d x_n - \sum_{i=1}^{n-1}y_i \d x_i,\]
with corresponding Liouville vector field
\[X_\rho= \frac{\rho(x_n)}{\rho'(x_n)} \del_{x_n} +
y_n \frac{\rho(x_n)\rho''(x_n)}{(\rho'(x_n))^2} \del_{y_n} + \sum_{i=1}^{n-1}y_i \del_{y_i},\]
that is gradient-like for the function $-x_n$. Near $x_n=1$, we have $\rho(x_n)=1-x_n$
and thus 
\[\lambda_\rho=(x_n-1)\d y_n - \sum_{i=1}^{n-1} y_i \d x_i.\]

Use coordinates on $\R^{2n-1}$ and $\S\R^{2n-1}$ as in Section \ref{sec:basics}
and consider the embedding germ $W_{\rho}\to \S \R^{2n-1}$ along $\{x_n=1\}$ given by
\begin{alignat*}{2}
&t=\log(1-x_n),\qquad  &&\zeta=-y_n,\\
&q_j=x_j,\qquad  &&p_j=\frac{y_j}{1-x_n},\qquad j=1,\dots,n-1.
\end{alignat*}
The pull-back of $\lambda_{\can}$ under this germ equals $\lambda_\rho$.
Since $X_{\rho}$ is gradient like for $-x_{n}$, the germ extends to an embedding
$\phi_\rho\colon W_\rho \to \S \R^{2n-1}$ with the property
$\phi_\rho^* \lambda_{\can}=\lambda_\rho$ and the extension is unique.
A straightforward calculation shows that the image of the embedding $\phi_\rho$
is $\{t<\log \rho(0)\}$.

Consider an exact Lagrangian filling $L \subset \S \R^{2n-1}$ of $\Sigma\subset\R^{2n-1}$.
We first push $L$ along the Liouville flow for sufficiently large negative time
to ensure that it is cylindrical in a neighborhood of the region $\{t \geq \ln\rho(0)\}$.
Let $h$ denote the compactly supported primitive of $\lambda_{\can}$ on $L$.
Then push $\phi_\rho^{-1}(L)$ along the flow
of $\del_z$ during the time $-h$. This gives a Legendrian submanifold
$\Lambda_\rho^+(L)\subset \R^{2n+1}$.

The diffeomorphism $\tau\colon\R^{2n+1} \to \R^{2n+1}$,
\[\tau(x_1,y_1,\dots,x_{n-1},y_{n-1},x_n,y_n,z)=(x_1,y_1,\dots,x_{n-1},y_{n-1},-x_n,-y_n,z) \]
satisfies $\tau^*\alpha = \alpha$. Thus $\tau(\Lambda_\rho^+(L))$ is Legendrian and we define
$\Lambda_\rho^-(L):=\tau(\Lambda_\rho^+(L))$. 

Since $\rho$ is even it follows that if $L$ and $L'$ are both exact Lagrangian fillings of $\Sigma$ 
then the union $\Lambda_\rho^-(L) \cup \Sigma \cup \Lambda_\rho^+(L')$ is a smooth Legendrian
submanifold diffeomorphic to the manifold resulting from gluing $L$ and $L'$ along $\Sigma$.
We write 
\[\Lambda_{\rho}(L,L'):=\Lambda_\rho^-(L) \cup \Sigma \cup \Lambda_\rho^+(L').\] 
If the function $\rho$ is deformed in the convex set of functions
satisfying conditions $(1)-(3)$ above, then the Legendrian $\Lambda_\rho(L,L')$ varies by
Legendrian isotopy. We therefore often drop the subscript $\rho$ from the notation and write
simply $\Lambda(L,L')$ when we discuss the Legendrian isotopy class.

\begin{rem}\label{r:frontspecial}
If we pick the function $\rho$ so that $\rho(u)=1-|u|$ in the region where $L$ is 
cylindrical, then for each branch $\{\zeta=h(q_1,\dots,q_{n-1})\}$ of the front projection
of $\Sigma$, the corresponding branch for $\psi(\Lambda_\rho^+(L))$ near $x_n=0$ satisfies
$z=\rho(x_n)h(x_1,\dots,x_{n-1})$. Hence the function $\rho$
scales the front projection of $\Sigma$ and the smoothing near $\{x_n=0\}$
can alternatively be described in terms of front projections. In particular,
the construction from \cite[Section 2.1]{ekholm_2016} is a special case
of the one presented above.
\end{rem}

\section{Proofs of Theorem \ref{t:main} and Corollary \ref{cor:unknot}}\label{sec:prooftmain}

We will use notation as in Sections \ref{sec:basics} and \ref{sec:glue}. Fix a function
$\rho$ which satisfies $(1)-(3)$ in Section \ref{sec:glue} and choose it to be
concave on $(0,1)$. 

\subsection{Proof of Theorem \ref{t:main}}

\subsubsection*{Step 1: A preliminary Legendrian isotopy.}
Note that the vector field $\del_\theta$ on $\R^{2n+1}$ is transverse to the hypersurface
$W_{\rho}\subset\R^{2n+1}$. For $s\in[0,1]$, consider the contact form on $\R^{2n+1}$ given by
\[\alpha_s=\frac{1}{1-s+sr^2/2}\left(\d z + \frac{r^2}{2} \d \theta\right)-\sum_{j=1}^{n-1}y_{j}\d x_{j},\] 
with corresponding Reeb vector field
\[R_{s}=(1-s)\del_z + s \del_\theta.\]
Let $\lambda_s=\alpha_s|_{W_{\rho}}$ be the restriction of $\lambda_{s}$ to $W_{\rho}$.
Flowing along the Liouville vector field $X_{\can}=\del_{t}$ of $\S\R^{2n-1}$ gives an
isotopy of embeddings $\psi_s\colon W_{\rho}\to \S \R^{2n-1}$ with the following properties:
\[ 
\psi_{0}=\phi_{\rho},\quad
\psi_s^*\lambda_{\can}=\lambda_s,\;0\le s\le 1,\quad\text{and }\quad
\psi_{1} \text{ is surjective}.
\]
Flowing $\psi_s(L)$ for time $-h$ along the Reeb vector field $R_{s}$ produces a Legendrian
isotopy $\Lambda_{s}^{+}(L)$, $0\le s\le 1$, that starts at $\Lambda_{0}^{+}(L)=\Lambda_{\rho}^{+}(L)$.
We have a directly analogous Legendrian isotopy that starts instead from $\Lambda_{s}^{-}(L)$.
In order for these isotopies to glue to an isotopy of $\Lambda_{\rho}(L,L)$ we must make sure
that the two halves do not intersect during the isotopy. It is easy to see that they do not
intersect if the function $h$ is sufficiently small since then $\Lambda^\pm_s(L)$ lies in
$\{\pm x_n>0\}$ for all $s\in[0,1]$. If we flow $L$ along the Liouville flow in negative time
then the function $h$ decreases exponentially. It follows that the desired isotopy exists.

\subsubsection*{Step 2: Construction of a pre-Lagrangian filling.}

\begin{figure}[!h]
\includegraphics[scale=0.7]{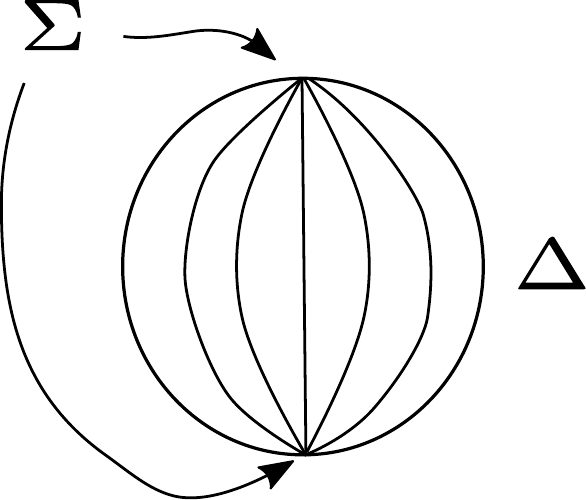}
\caption{The pre-Lagrangian disk $\Delta$ with its singular Legendrian foliation. The singular set $\Sigma$
is an equatorial sphere in $\del \Delta$.}
\label{fig:prelag}
\end{figure}
After this preliminary isotopy, we observe that $\Lambda(L,L) \subset \R^{2n+1}$ bounds
a pre-Lagrangian submanifold (see \cite{ehs} for this notion) $\Delta$ which is foliated by $\varphi_{\del_{\theta}}^t(\Lambda^+(L))$ for
$t \in [0,\pi]$, see Figure \ref{fig:prelag}. The singular set of the characteristic foliation of $\Delta$
is precisely $\Sigma=\Delta \cap \{r=0\}$.
\subsubsection*{Step 3: Isotoping pre-Lagrangian fillings}
Any two Lagrangian embeddings (with Legendrian boundary) of $L$ lead
to pre-Lagrangian fillings $\Delta$ with diffeomorphic Legendrian foliations.
We show that they are contact isotopic. 

Consider two such embeddings $f,g\colon\Delta\to \R^{2n+1}$. By assumption,
the induced trivializations of $T L \otimes \C$ are homotopic and Gromov's h-principle
for Legendrian immersions implies that the restrictions $f|_{L_{0}}$ and $g|_{L_{0}}$ to the central
leaf $L_0$ of $\Delta$ are Legendrian regularly homotopic. By general position, the restriction of this regular homotopy
to the codimension $1$ skeleton of $L_0$ is an isotropic isotopy. Retracting $L_0$ to its codimension $1$ skeleton, we then find that $f|_{L_0}$ and $g|_{L_0}$ are Legendrian isotopic. Extending this Legendrian isotopy to an ambient contact isotopy, we reduce to the case that $f|_{L_{0}}=g|_{L_{0}}$.

Assume then that $f|_{L_{0}}=g|_{L_{0}}$. We first deform $g$ slightly near the singular set $\Sigma$ of the Legendrian foliation of $\Delta$ so that $f$ and $g$ coincides there as follows. Pick a germ
of Legendrian foliation that contains $L_0$ as a leaf and is transverse to
all leaves of $f(\Delta)$ and of $g(\Delta)$ near $\Sigma$. The corresponding front projections of these leaves near $\Sigma$ then consist of families of graphs of functions (all vanishing to first order on $\Sigma$) and a linear interpolation between these families give the desired isotopy near $\Sigma$.

We consider next a $1$-jet space neighborhood of $L_0$, here nearby leaves correspond to a family of $1$-jet graphs of functions $h_t\colon L_0\to \R$ which are increasing with respect to $t$ and coincide near $\del L_0$.
Linear interpolation connects any two such families of functions within the class of
families that are increasing in $t$. Using the corresponding isotopies, we may then assume
that $f$ and $g$ agree near the central leaf $L_0$.

To finish the proof, identify $\Delta \setminus \Sigma$ with $L_0\times [-\frac{\pi}{2},\frac{\pi}{2}]$ with the
foliation given by $L_0\times\{\theta\}$, $\theta \in [-\frac{\pi}{2},\frac{\pi}{2}],$ and consider the isotopy
$\phi_s\colon\Delta\setminus \Sigma \to \Delta\setminus \Sigma$ given by $\phi_s(x,\theta)=(x,(1-s+s\epsilon)\theta)$.
For $\epsilon >0$ small enough, $f$ and $g$ agree on $\phi_1(\Delta\setminus \Sigma)$.
The isotopy $f\circ \phi_s$ extends to a contact isotopy $F_s$ since it
has a fixed characteristic foliation. Similarly, $g\circ \phi_s$ extends to a contact isotopy $G_s$.
Moreover, we may assume that $F_s=G_s$ near $f(\Sigma)=g(\Sigma)$ since $f$ and $g$
agree there. Then $F_s^{-1}\circ G_s \circ g$ equals $g$ for $s=0$ and $f$ for $s=1$,
and $F_s^{-1}\circ G_s$ extends smoothly to $f(\Sigma)$ by the identity.
\qed

\subsection{Proof of Corollary \ref{cor:unknot}}
Since all trivializations of $T D \otimes \C$ are homotopic, Theorem \ref{t:main} shows
that the Legendrian isotopy class of $\Lambda(D,D)$ is independent of $D$.
Consider the case when $D$ is a trivial disk filling. Using the construction in \cite[Section 2.1]{ekholm_2016}, see Remark \ref{r:frontspecial}, we represent $D$ as half of the front of the standard unknot and then see directly that $\Lambda(D,D)$ is an unknot.\qed

\section{Loose Legendrian spheres}\label{sec:loose}
\begin{figure}[h!]
	\includegraphics[width=0.3\linewidth]{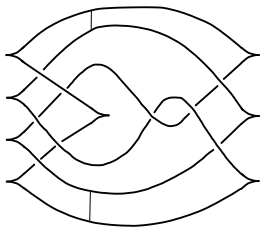}
	\caption{Front projection of the Legendrian $\overline{9_{46}}$. The two places to pinch to get the distinct disk fillings $D_0$ and $D_{1}$ are indicated as vertical lines.}
	\label{fig:knot}
\end{figure}

In this section we present Emmy Murphy's argument showing that $\Lambda(D_0,D_1)$ is loose for distinct disk fillings of the Legendrian knot in the topological isotopy class $\overline{9_{46}}$. The knot is shown in Figure \ref{fig:knot}. The two distinct disk fillings of the knot are obtained by pinching the two Reeb chords indicated in Figure \ref{fig:knot}. 

\begin{figure}[htb]
	\includegraphics[width=.45\linewidth]{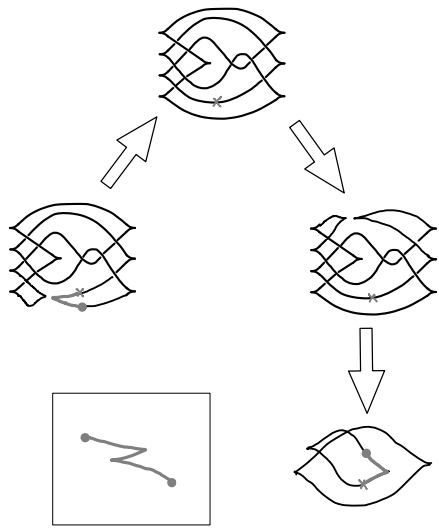}
	\caption{The embedded zig-zag in the front of $\Lambda(D_0,D_1)$.  }
	\label{fig:loose}
\end{figure}

Figure \ref{fig:loose} shows a zig-zag path, embedded in the front. More precisely, the leftmost front diagram shows a slice in the front of $\Lambda(D_0,D_1)$ which lies in $D_0$ after the pinch. The path starts at the gray dot indicated there and moves up the left cusp to the cross. The path then continues as the cross in slices of $D_0$ approaching the central slice shown at the top of the figure, and down in the slices corresponding to $D_1$ until it passes the pinching slice, as indicated in the next front. Then the slices of $D_1$ undergo isotopy and ends in the lower right front. Here the cross goes up the right cusp completing the zig-zag. After deformation of the front along the zig-zag path we find an embedded zig-zag as shown in the lower left of the figure which then gives an explicit loose chart. 

\bibliographystyle{abbrv}
\bibliography{biblio}
\end{document}